\documentclass[12pt]{article}
\usepackage[margin=1.25in]{geometry}

\usepackage[utf8]{inputenc}
\usepackage[T1]{fontenc}
\usepackage{lmodern}
\usepackage{etex}
\reserveinserts{28}
\usepackage{authblk}
\usepackage{amsthm}
\usepackage{amssymb}
\usepackage{amsmath}
\usepackage{float}
\usepackage{tikz}
\usetikzlibrary{arrows,calc,shapes,decorations.pathreplacing}
\usepackage{hyperref}
\hypersetup{pdftex,colorlinks=true,allcolors=blue}
\usepackage{soul}

\newtheorem{theorem}{Theorem}
\newtheorem{lemma}[theorem]{Lemma}

\newtheorem{remark}{Remark}

\def\W{\mathrm{W}}

\newcommand{\eq}[1]{(\ref{eq:#1})}

\begin{document}
\title{On reflected L\'evy processes with collapse}

\author{Onno Boxma\thanks{EURANDOM and Department of Mathematics and Computer Science, Eindhoven University of Technology, the Netherlands. \texttt{o.j.boxma@tue.nl}}
\thanks{Partly funded by the NWO Gravitation project N{\sc etworks}, grant number 024.002.003.} ,
Offer Kella\thanks{Department of Statistics, The Hebrew University of Jerusalem; Jerusalem 9190501, Israel. \texttt{offer.kella@gmail.com}} 
\thanks{Partially supported by ISF grant 3336/24 and  the Vigevani Chair in Statistics.} \ and David Perry\thanks{Industrial Engineering and Technology Management, Holon Institute of Technology, Holon 5810201, Israel. \texttt{davidper@hit.ac.il}} \thanks{Partially supported by ISF grant 3336/24.}
}

\date{\today}

\maketitle

\begin{abstract}
We consider a L\'evy process reflected at the origin with additional i.i.d.\ {\em collapses} that occur at Poisson epochs, where a {\em collapse} is a jump downward to a state which is a random fraction of the state just before the jump. We first study the general case, then specialize to the case where the L\'evy process is spectrally positive and finally we specialize further to the two cases where the L\'evy process is a Brownian motion and a compound Poisson process with exponential jumps minus a linear slope.
\end{abstract}


\bigskip
\noindent {\bf AMS Subject Classification (MSC2020):} Primary 60G51; Secondary 60K25.

\section{Introduction}

This paper is devoted to the study of L\'evy processes which are reflected at the origin, and which have the special feature that at Poisson epochs they are subject to random collapses -- jumps downward where the jump size is a random proportion of the state just before the jump. Such processes occur naturally in several different fields. To name a few:
\begin{itemize}
    \item {\bf Population genetics.} Population sizes may fluctuate according to birth-and-death processes, and experience disasters at random epochs. After appropriate scaling such a process may be approximated by a reflected Brownian motion subject to collapse. This case is given special attention in Example~\ref{sec6.1} of the current study.
    \item {\bf Geophysics.} Pressure in earth layers may fluctuate randomly in between earthquakes or volcano eruptions (collapses).
    \item {\bf TCP (Transmission Control Protocol).}
    TCP adapts the window size (transmission rate) of data transfers to the congestion in the network. The window size along a path is increased until a signal is received that the path gets too congested; the window size is then proportionally reduced (additive increase multiplicative decrease). 
    Stochastic processes with random collapses have been intensively studied in the analysis of the Transmission Control Protocol of the Internet (see, among many others, \cite{AABD,DGR2002,GRZ2004,LvL}), but typically those processes are non-decreasing (often linear) in between collapses.
    \item {\bf Perishable inventories.}  In stochastic models of perishable inventory systems (PIS) with random input, the inventory typically increases due to production and decreases because of demand. However, in the sizable literature of such PIS (cf.\ the survey \cite{BPSreview}) not much attention has been given to the modeling of incidents like power outages, which may make a random proportion of the stored goods (e.g., food or blood) obsolete.  See \cite{KPS} for the case of PIS with {\em clearing}, i.e., collapses in which all the goods become obsolete. 
    \item {\bf Cash management}. We dwell a bit longer on this particular application area, as it may be less known than the ones above, while it formed one of our main motivations for the present study.
    The \textit{Peer-to-Peer} (P2P) lending platform is the practice of
lending money through online services that match lenders with borrowers. The
P2P firm offers a platform where individuals or businesses can lend
directly to other individuals or businesses without the need for a bank as a
middleman. The firm takes brokerage fees for providing this match-making
platform. Compared to investment and savings products from banks, borrowers
can borrow money at lower interest rates, and lenders can earn higher
returns. Although the P2P lending firm applies a strict screening system, it
may happen that a borrower is in arrears. If that happens, the P2P firm
initiates collection procedures against the debtor. P2P systems are exposed
only to the exceptional volatility caused by a severe crash in which many
borrowers simultaneously go bankrupt. Historical examples of such crashes
are the dot.com bubble in 2000, the sub-prime crisis in 2008 that started
with the crash of Lehman Brothers, the corona crisis of 2020, etc. In these
cases, the stock market crash was accompanied by bankruptcies in the entire
business sector of small and medium businesses and even large companies. In
the events mentioned above the crisis hit the entire capital market, since
all those economic entities could not meet 
their financial obligation.
In \cite{BPSP2P} a stochastic analysis of such cash management processes with crashes/collapses is presented, with a compound Poisson process in between collapses. However, the cash management process generated by big businesses together with many "small" individuals can be more naturally approximated by an independent sum of a compound Poisson and a Brownian component, reflected at the origin, which is a special case of the more general reflected L\'evy process considered in this paper.
\end{itemize}

For some background on growth collapse processes, where, unlike in the current paper, the processes are nondecreasing between collapse epochs, one may sample, among others, \cite{BKP2011,BPSZ2006,HKPS2023,K2009,KL2010,LS2011} and references therein.

Globally, our {\em goals} are to study the stability conditions of reflected L\'evy processes with random collapses, and to determine the stationary distribution (if it exists) when the L\'evy process is spectrally positive (i.e., the L\'evy process itself has no negative jumps).

Our {\em main contributions} are:
\begin{itemize}
    \item 
    A careful discussion of the conditions for the existence of a limiting/ergodic/
    \\
    stationary distribution of the reflected L\'evy process with collapses.
\item
For the case of a spectrally positive L\'evy process, we determine the LST (Laplace-Stieltjes Transform) and moments of that limiting distribution.
\item
For the special case of a compound Poisson process with regularly varying jump sizes, we prove that the tail of the limiting distribution is just as heavy as the tail of the jump size distribution.
\item
For the special cases of Brownian motion and of a compound Poisson process with exponentially distributed jump sizes (corresponding to an $M/M/1$ queue with collapses) the LST of the limiting distribution is worked out in more detail, revealing a relation to incomplete Beta functions.
\end{itemize}

\begin{remark}{\rm
    There are many extensions of our model possible; we mention a few. (i)~The intervals between collapses could have a general distribution (cf.\ \cite{BKP2011}, where the resulting model is studied for the case of a compound Poisson process minus a drift).
    (ii)~It would be interesting to study reflected L\'evy processes in which work decreases proportionally ("shot noise");
    cf.\ \cite{K2009}, where the relation between collapses and shot noise is explored. (iii)~In perishable inventory systems, it is often assumed that items have a finite life time \cite{BPSreview}. In our setting, this could give rise to {\em doubly-reflected} L\'evy processes with collapses. (iv)~Next to proportional collapses, one might also allow state dependent jumps {\em upward}. In cash management and the crypto exchange market, such jumps might occur, although they seem milder than the collapses.
    }
\end{remark}

The remainder of the paper is organized as follows. In Section~\ref{secmodel} we introduce the model under consideration: the L\'evy process, the collapse procedure that is applied to it, the reflected version of the resulting process, and its local time at zero. 
We also show that the process level just before successive collapses satisfies a Lindley-type recursion of an autoregressive sequence.
Section~\ref{secLindley} considers this recursion in detail. It leads to the conclusion that, if a limiting distribution for this sequence exists, then the limiting random variable can be decomposed into two independent terms.
In Section~\ref{secstability} we discuss the conditions under which the reflected L\'evy process with uniform collapses at Poisson epochs has a limiting/ergodic/stationary distribution. By PASTA, this is also the limiting distribution of the process level just before collapses.
Section~\ref{secspectral} presents the derivation of this limiting distribution.
The special cases of reflected Brownian motion with collapses, of an $M/M/1$ queue with collapses, and of a heavy-tailed $M/G/1$ queue with collapses are discussed in Section~\ref{secexample}.

\section{The general model}
\label{secmodel}
In what follows $x\vee y=\max(x,y)$, $x\wedge y=\min(x,y)$, $x^+=x\vee 0$ and $x^-=(-x)^+$. Also, a.s., w.l.o.g., iff and PASTA, abbreviate: almost surely, without loss of generality, if and only if and Poisson Arrivals See Time Averages, respectively. Moreover $\sim$ denotes "distributed" or "distributed like" and $\approx$ denotes "asymptotic to".

Consider a c\`adl\`ag L\'evy process $X=\{X_t|\,t\ge 0\}$ with $X_0=0$. Associated with each L\'evy process is a triplet $(c,\sigma^2,\nu)$ where $c\in\mathbb{R}$, $\sigma\ge 0$ and $\nu$ is a (sigma-finite) measure, called the {\em L\'evy measure} which satisfies $\nu(\{0\})=0$, $\int_{\mathbb{R}}\,x^2\wedge 1\,\nu(dx)<\infty$. Any such L\'evy process can be decomposed into a sum of independent processes, where one is a Brownian motion (possibly only a drift), another is a compound Poisson process (possibly zero) with jumps in $\mathbb{R}\setminus[-1,1]$ and finally a convergent sum of centered compound Poisson processes of the form $\sum_{n=1}^\infty (X_{n,t}-tEX_{n,1})$, where the absolute values of the jumps of the process $X_n$ are in $\left(\frac{1}{n+1},\frac{1}{n}\right]$. For background on L\'evy processes and L\'evy queues see, e.g., \cite{Bertoin,DM,KPS, Kyprianou,Sato}.

In addition, let $\tau,\tau_1,\tau_2,\ldots$ be i.i.d.\ $\exp(\lambda)$ distributed random variables which are also independent of $X$ and denote $S_0=0$ and $S_n=\sum_{i=1}^n\tau_i$ for $i\ge 1$, with renewal counting (Poisson) process $N_t=\sup\{n|\,S_n\le t\}$. Finally, let $U,U_1,U_2,\ldots$ be i.i.d.\ with $P(U\in[0,1])=1$. These are thought of as random proportions and are independent of all the other random objects defined thus far.

For $t\ge 0$, denote $L_t=-\inf_{0\le s\le t}X_s$, $W_t=X_t+L_t$ and, for $x\ge 0$,
\begin{align}\label{eq:reflection}
L^x_t&=\left(-\inf_{0\le s\le t}(x+X_s)\right)^+=(L_t-x)^+ , \nonumber\\
W^x_t&=x+X_t+L^x_t=X_t+x\vee L_t=W_t+ (x-L_t)^+\,.
\end{align}
This is the reflection (Skorokhod) map associated with $\{x+X_t|\, t\ge 0\}$ where $(L^x_t,W^x_t)$ are the unique c\`adl\`ag processes that jointly satisfy (see \cite{K2006}):
\begin{enumerate}
\item $L^x_0=0$ and $L^x_t$ is nondecreasing in $t$,
\item $W^x_t\ge 0$ for each $t\ge 0$ and
\item $\{t|\,W^x_t=0\}\subset \{t|\,L^x_s<L^x_t\ \forall s\in [0,t)\}$.
\end{enumerate}

Now denote $\sigma_sX=\{X_{t+s}-X_s|\,t\ge 0\}$ (which is independent of the history of $X$ until time $s$) and let $(\sigma_sL^x_t,\sigma_sW^x_t)$ be the corresponding reflection map associated with $\sigma_sX$.
Let $Z=\{Z_t|\,t\ge 0\}$ be defined as follows. Assume that we have already defined it on $[0,S_n)$. Then for $t\in [S_n,S_{n+1})$, $Z_t=\sigma_{S_n}W^{U_nZ_{S_n-}}_{t-S_n}$, where
$Z_{t-}=\lim_{v\uparrow t}Z_v$. Namely, at time $S_n$, which we refer to as {\em collapse epochs}, the content modeled by the process $Z$ collapses to a fraction $U_n$ of its precollapse level and continues according to a reflected L\'evy process until the next collapse epoch. $Z_0$ may have any distribution and is assumed to be independent of everything else. In fact, because of the assumptions made until now, $Z$ is a (nonnegative) Markov process with the following generator:
\begin{align}
\mathcal{A}g(x)&=cg'(x)+\frac{\sigma^2}{2}g''(x)+\int_{\mathbb{R}}\left(g(x+y)-g(x)-g'(x)y1_{[-1,1]}(y)\right)\nu(dy)\nonumber\\
&+\lambda(Eg(Ux)-g(x))\,
\end{align}
for twice continuously differentiable $g$ with $g(x)=g(0)$ for $x\le 0$ (hence, $g'(x)=g''(x)=0$ for $x\le 0$). In particular, for $x\le 0$ it may be checked that
\begin{equation}
\mathcal{A}g(x)=\int_{(-x,\infty)}(g(x+y)-g(0))\nu(dy)=\int_0^\infty g'(y)\nu(-x+y)\,dy\,.
\end{equation} This can be shown by applying the (local) martingales discussed in \cite{KY2013} but will not be needed for what follows.

A famous result, which in particular may be found in \cite{Bertoin,DM,Kyprianou,Sato}, is that $W_\tau$ and $L_\tau$ are independent. By \eq{reflection} this implies the following decomposition:
\begin{equation}
(W^x_\tau,L^x_\tau)=(W_\tau,0)+((x-L_\tau)^+,(L_\tau-x)^+) ,
\end{equation}
where the two random vectors on the right are independent. In particular, denoting $w(\alpha)=Ee^{-\alpha W_\tau}$, we have that, for $\alpha,\beta\ge 0$,
\begin{equation}\label{eq:decomposition}
Ee^{-\alpha W^x_\tau-\beta L^x_\tau}=w(\alpha)\,Ee^{-\alpha (x-L_\tau)^+-\beta (L_\tau-x)^+}\,.
\end{equation}

With our setup, it follows from \eq{decomposition} that
\begin{equation}
Z_{S_n-}=\sigma_{S_{n-1}} W_{\tau_n}+(Z_{S_{n-1}-}U_n-\sigma_{S_{n-1}}L_{\tau_n})^+ ,
\end{equation}
where $\sigma_{S_{n-1}} W_{\tau_n},Z_{S_{n-1}-},U_n,\sigma_{S_{n-1}}L_{\tau_n}$ are independent and $\sigma_{S_{n-1}} W_{\tau_n},U_n,\sigma_{S_{n-1}}L_{\tau_n}$ are distributed like $W_\tau,U,L_\tau$, respectively.

We observe that by PASTA the ergodic distribution of the continuous time process $Z$ is the same as that of the discrete time process $Z_{S_n-}$ and therefore it is important to study this latter process. With $\zeta_n=Z_{S_n-}$, $V_n=\sigma_{S_{n-1}}W_{\tau_n}$ and $Y_n=\sigma_{S_{n-1}}L_{\tau_n}$ we therefore have the following Lindley-style version of an autoregressive sequence:
\begin{equation}
\zeta_n=V_n+(\zeta_{n-1}U_n-Y_n)^+\,,
\end{equation}
where all four random variables on the right are independent with $V_n\sim W_\tau$, $U_n\sim U$ and $Y_n\sim L_\tau$.
In the next section we consider this recursion in some detail.

\section{Lindley-style autoregressive recursions}
\label{secLindley}
Consider the (for now, deterministic) recursion 
\begin{equation}\label{eq:linaut}
z_n=v_n+(z_{n-1}u_n-y_n)^+ ,  ~~~~ n \geq 1,
\end{equation}
where $v_{i-1},y_i\ge 0$, $u_i\in[0,1]$ for $i\ge 1$. Also, let $z_0=v_0$. This is the same as letting $z_{-1}=0$, defining $u_0\in[0,1]$ arbitrarily  and starting the recursion one index earlier. Also let $\Pi_0=1$ and $\Pi_n=\prod_{i=1}^n u_i$ for $n\ge 1$. 

We first observe that if $z'_n=v_n+(z_{n-1}'u_n-y_n)^+$ and $z''_n=v_n+(z_{n-1}''u_n-y_n)^+$, for any two choices where $z''_0\ge z'_0\ge 0$, then it follows by induction that $z''_n\ge z'_n$ for all $n\ge 0$. Moreover: 
\begin{itemize}
\item When $z''_{n-1}u_n\le y_n$ we have that $z''_n-z'_n=0\le (z''_{n-1}-z'_{n-1})u_n$.
\item When $z'_{n-1}u_n>y_n$ we have that $z'_n-z_n=(z''_{n-1}-z'_{n-1})u_n$.
\item When $z'_{n-1}u_n\le y_n<z''_{n-1}u_n$ we have that
\begin{equation}
z''_n-z'_n=z''_{n-1}u_n-y_n\le z''_{n-1}u_n-y_n+y_n-z'_{n-1}u_n=(z''_{n-1}-z'_{n-1})u_n\,.
\end{equation}
\end{itemize}
Therefore, for all $n\ge 1$, $z''_n-z'_n\le (z''_{n-1}-z'_{n-1})u_n$ and thus, for every $n\ge 0$,
\begin{equation}\label{eq:insensitive}
 0\le z''_n-z'_n\le (z''_0-z'_0)\Pi_n\,.
\end{equation}
In particular this implies:
\begin{lemma}
If $\Pi_n\to0$ as $n\to\infty$ then $z'_n-z_n\to 0$ as $n\to\infty$ for any choice of $z'_0\ge 0$ (smaller, equal or larger than $z_0=v_0$).
\end{lemma}
We note that since we will be replacing $u_i$ by i.i.d.\ random variables $U_i\sim U$ with support in $[0,1]$, then, unless $P(U=1)=1$, it will always hold that $\prod_{i=1}^nU_i\to 0$ as $n\to\infty$.

We begin by assuming that $u_i>0$ for all $i\ge 1$, so that $\Pi_n>0$ for all $n\ge1$.
Denoting, for $n\ge 1$,
\begin{equation}
w_n=\frac{z_n-v_n}{\Pi_n}, \qquad x_n=\frac{v_{n-1}}{\Pi_{n-1}}-\frac{y_n}{\Pi_n},
\end{equation}
we have that $w_0=0$ and, for $n\ge 1$, $w_n=(w_{n-1}+x_n)^+$, which is Lindley's recursion. Thus, with $s_0=0$ and $s_n=\sum_{i=1}^nx_i$ for $n\ge 1$ we have the well known solution:
\begin{equation}
w_n=s_n-\min_{0\le k\le n}s_k=\max_{0\le k\le n}(s_n-s_{n-k})
\end{equation}
and thus
\begin{equation}
z_n=v_n+\max_{0\le k\le n} ((s_n-s_{n-k})\Pi_n)\,.
\end{equation}
For $k=0$ we have that $(s_n-s_{n-k})\Pi_n=0$ while for $1\le k\le n$ we have
\begin{align}
(s_n-s_{n-k})\Pi_n&=\sum_{i=n-k+1}^n \left(\frac{v_{i-1}}{\Pi_{i-1}}-\frac{y_i}{\Pi_i}\right)\Pi_n
\nonumber\\
&=\sum_{i=1}^k \left(\frac{v_{n-i}}{\Pi_{n-i}}-\frac{y_{n-i+1}}{\Pi_{n-i+1}}\right)\Pi_n\,,
\end{align}
recalling that $v_0=0$ (relevant for $i=k=n$). Now,
\begin{equation}
\frac{\Pi_n}{\Pi_{n-i}}=\prod_{j=n-i+1}^nu_j=\prod_{j=1}^i u_{n-j+1}
\end{equation}
and similarly, for $i=1$, we have that $\Pi_n/\Pi_{n-1+1}=1$ and for $2\le i\le k\le n$,
\begin{equation}
\frac{\Pi_n}{\Pi_{n-i+1}}=\prod_{j=n-i+2}^n u_j=\prod_{j=1}^{i-1} u_{n-j+1}\,.
\end{equation}
This implies that we can write,
\begin{equation}
(s_n-s_{n-k})\Pi_n=\sum_{i=1}^k v_{n-i}\prod_{j=1}^iu_{n-j+1}-\sum_{i=1}^ky_{n-i+1}\prod_{j=1}^{i-1}u_{n-j+1}\,,
\end{equation}
where an empty product is defined to be one and an empty sum is zero. Therefore, we finally have that
\begin{equation}\label{eq:explicit}
z_n=v_n+\max_{0\le k\le n}\left(\sum_{i=1}^k(v_{n-i}u_{n-i+1}-y_{n-i+1})\prod_{j=1}^{i-1}u_{n-j+1}\right)\,.
\end{equation}
We emphasize the fact that for $n=0$ we have that (since empty sums are zero), \eq{explicit} implies that $z_0=v_0$ as initially assumed.
Finally, it is a straightforward exercise to show by induction that \eq{explicit} holds also for the case where $u_i$ are allowed to be zero. For each such $i$ we have that $z_i=v_i$ and the recursion evidently restarts from that point.

Now, if we replace $\{(v_{i-1},u_i,y_i)|\,i\ge 1\}$ by three independent i.i.d.\ sequences $\{(V_{i-1},U_i,Y_i)|\,i\ge 1\}$ of nonnegative random variables and $\{z_n|\,n\ge 0\}$ by the random variable $\zeta_n$, then by taking each of the random vectors $V_0,\ldots,V_n$, $U_1,\ldots,U_n$ and $Y_1,\ldots,Y_n$ in reverse order it follows that
\begin{equation}\label{eq:zetan}
\zeta_n\sim V_0+\max_{0\le k\le n}\left(\sum_{i=1}^k (V_iU_i-Y_i)\prod_{j=1}^{i-1}U_j\right)\,,
\end{equation}
where an empty sum is zero and and empty product is one.
Therefore, $\zeta_n$ is stochastically increasing and bounded above by $\sum_{i=0}^nV_i\prod_{j=1}^iU_j$ and hence, if $EV_0<\infty$ and $P(U_i=1)<1$, $\zeta_n$ converges in distribution to an a.s.\ finite random variable $\zeta$ which is distributed like an independent sum of two random variables.
The first is distributed like $V_0$ and the second like
\begin{equation}\label{eq:Wpart}
\max_{k\ge 0}\left(\sum_{i=1}^k (V_iU_i- Y_i)\prod_{j=1}^{i-1}U_j\right)\,.
\end{equation}
This is a generalization of the setup (hence the conclusions) considered in \cite{BKM2023}. See also \cite{BLMP2021} for a discussion of recursions of the form $W_n=(W_{n-1}U_n+X_n)^+$.

\section{Stability of the continuous time process}
\label{secstability}
An important question is what are the conditions on $Z$ that ensure that there exists a limiting/ergodic/stationary distribution. If $P(U=1)=1$ then we have a reflected L\'evy process, which is well understood and has been discussed many times in the literature and textbooks, and thus we will assume that $P(U=1)<1$. Since the collapse epochs occur according to a Poisson process with rate $\lambda$, it is clear that true collapses (with $U_i<1$) occur according to a Poisson process with rate $\lambda P(U<1)$. Therefore, w.l.o.g., let us assume that $P(U=1)=0$. If $P(U=0)>0$, then after a geometrically distributed number of trials the process jumps to zero (often called: {\em clearing}), which implies that it is regenerative and therefore a proper limiting=ergodic distribution clearly exists. It remains to consider the case $P(U=0)=P(U=1)=0$. In this case it is clear that $\prod_{i=1}^\infty U_i=0$ a.s.

It is easy to check that if $x\le y$ then $W^x_t\le W^y_t$ for all $t\ge 0$ and that $W^y_t-W^x_t$ is nonincreasing in $t$. See \cite{KW1996} for results like this in the multidimensional setting.
Therefore, setting $Z^x$ to be the process $Z$ when started from $Z_0=x$ we have that
\begin{equation}
Z_{\tau_1}^y-Z_{\tau_1}^x=(W_{\tau_1-}^y-W_{\tau_1-}^x)U_1\le (y-x)U_1\,,
\end{equation}
where the left hand side is nonnegative. Hence by induction it also follows that
\begin{equation}
0\le Z_{S_n}^y-Z_{S_n}^x\le (y-x)\prod_{i=1}^n U_i
\end{equation}
where the right hand side vanishes a.s. as $n\to\infty$. Therefore we also have that
\begin{equation}
0\le Z^y_t-Z^x_t\le (y-x)\prod_{i=1}^{N_t}U_i
\end{equation}
which vanishes a.s. as $t\to\infty$. Therefore, it follows that if a limiting distribution exists, it does not depend on the initial distribution of $Z_0$.

Since $W^x_t=X_t+L_t\vee x\le x+X_t+L_t=x+W_t$ it follows that $Z_{S_n}$ is stochastically bounded by $x\prod_{i=1}^nU_i+\sum_{i=1}^n V_i\prod_{j=i}^nU_j$ where all the variables on the right are independent and $V_i\sim W_\tau$. By results from \cite{K2009} this implies that $\{Z_{S_n}|\, n\ge 0\}$ and hence $\{Z_t|\,t\ge 0\}$ are bounded above by a process that has a stationary version. Therefore, tightness of these processes is assured. It also implies that $\{Z_t|\,t\ge 0\}$ has a stationary version by following the (Loynes-type) construction from \cite{K2009} which is there applied to the case where $W$ is nondecreasing.

Therefore we can apply PASTA to conclude that the (unique) ergodic distribution of $Z$ is the same as that of the discrete time process $\zeta_n=Z_{S_n-}$, where we recall $\zeta_n,\zeta$ from Sections~\ref{secmodel} and \ref{secLindley}. With $Z^*=\zeta$ having the steady state distribution of $\zeta_n$ (and hence, also of the process $Z$), we have that 
\begin{equation}\label{eq:wzutau}
Z^*\sim W^{Z^*U}_\tau\sim W_\tau+(Z^*U-L_\tau)^+\,,
\end{equation}
where $Z^*,W_\tau,U,L_\tau$ on the right are independent.

\section{The case where $X$ is spectrally positive}
\label{secspectral}

When $X$ is spectrally positive, that is, $\nu(-\infty,0]=0$, then with
\begin{equation}
\varphi(\alpha)=\log Ee^{-\alpha X_1}=-c\alpha+\frac{\sigma^2}{2}\alpha^2+\int_{(0,\infty)}\left(e^{-\alpha x}-1+\alpha x1_{(0,1]}(x)\right)\,\nu(dx) ,
\label{Laplaceexp}
\end{equation}
it may be concluded (e.g., applying the martingale from \cite{KW1992}, as is done in Thm.~3.10, p. 259 of \cite{A2003}) that if $X$ is not a subordinator and $\alpha_\lambda$ is the unique positive root of $\varphi(\alpha)=\lambda$, then $L_\tau\sim \exp(\alpha_\lambda)$ and 
\begin{equation}\label{eq:LSTwtau}
w(\alpha)=Ee^{-\alpha W_\tau}=\frac{1-\frac{\alpha}{\alpha_\lambda}}{1-\frac{\varphi(\alpha)}{\lambda}}\,.
\end{equation}
Therefore, it follows either by Thm.~3.10, p.~259 of \cite{A2003} or from \eq{decomposition} here that
\begin{equation}\label{eq:xexp}
Ee^{-\alpha W^x_\tau-\beta L^x_\tau}=w(\alpha)\,\frac{e^{-\alpha x}-\frac{\alpha+\beta}{\alpha_\lambda+\beta}e^{-\alpha_\lambda x}}{1-\frac{\alpha}{\alpha_\lambda}}\,,
\end{equation}
while by l'H\^opital's rule and continuity, for $\alpha=\alpha_\lambda$ the right hand side is
\begin{equation}
\frac{x+\frac{1}{\alpha_\lambda+\beta}}{\varphi'(\alpha_\lambda)}\,\lambda e^{-\alpha_\lambda x} .
\end{equation}
In particular, for $\alpha\not=\alpha_\lambda$,
\begin{equation}\label{eq:Wx}
Ee^{-\alpha W^x_\tau}=w(\alpha)\,\frac{e^{-\alpha x}-\frac{\alpha}{\alpha_\lambda}e^{-\alpha_\lambda x}}{1-\frac{\alpha}{\alpha_\lambda}}=\frac{e^{-\alpha x}-\frac{\alpha}{\alpha_\lambda}e^{-\alpha_\lambda x}}{1-\frac{\varphi(\alpha)}{\lambda}}\,,
\end{equation}
where the term to the right of $w(\alpha)$ is the LST of $(x-L_\tau)^+$. The case where $X$ is a subordinator is analyzed in \cite{K2009,BKP2011} and is thus ignored here. 

From \eq{wzutau} and \eq{Wx} it therefore follows that if we set $f(\alpha)=Ee^{-\alpha Z^*}$ then
\begin{equation}\label{eq:falpha}
f(\alpha)=\frac{Ef(\alpha U)-\frac{\alpha}{\alpha_\lambda} Ef(\alpha_\lambda U)}{1-\frac{\varphi(\alpha)}{\lambda}}\,.
\end{equation}

Noting that for $n\ge 1$ we have, with $\delta_{1n}=0$ for $n\ge 2$ and $\delta_{11}=1$, that
\begin{equation}\label{eq:EZstar0}
    \left(1-\frac{\varphi(\alpha)}{\lambda}\right)f^{(n)}(\alpha)-\frac{1}{\lambda}\,\sum_{k=0}^{n-1}{n\choose k}f^{(k)}(\alpha)\varphi^{(n-k)}(\alpha)=EU^nf^{(n)}(\alpha U)-\frac{Ef(\alpha_\lambda U)}{\alpha_\lambda}\delta_{1n}
\end{equation}
from which the following recursion is immediate, once we argue that $E(Z^*)^n<\infty$ whenever $\varphi^{(n)}(0)$ is finite.
\begin{lemma}\label{lem:moments}
    Assume that $P(U=1)<1$ and for $n\ge 0$ denote $m_n=E(Z^*)^n=(-1)^nf^{(n)}(0)$, and
    \begin{equation}
c_n=(-1)^n\varphi^{(n)}(0)=        \begin{cases}
            0&n=0\,\\
            c+\int_{(1,\infty)}x\nu(dx)&n=1\,\\
            \sigma^2+\int_{(0,\infty)}x^2\nu(dx)&n=2\,\\
            \int_{(0,\infty)}x^n\nu(dx)&n\ge 3\,.
        \end{cases}
    \end{equation}
    Then, for all $n\ge 0$, $m_n<\infty$ when $c_n<\infty$ and for $n\ge 2$
    \begin{equation}\label{eq:EZstarn}
        m_n=\frac{1}{1-EU^n}\,\frac{1}{\lambda}\,\sum_{k=0}^{n-1}{n\choose k}m_kc_{n-k},
    \end{equation}
    where $m_0=1$ and
    \begin{equation}\label{eq:EZstar1}
        m_1=\frac{1}{1-EU}\left(\frac{c_1}{\lambda}+\frac{Ef(\alpha_\lambda U)}{\alpha_\lambda}\right).
    \end{equation}
\end{lemma}

\medskip
\begin{proof}
    Since \eq{EZstarn} and \eq{EZstar1} are obtained by letting $\alpha\downarrow 0$ in \eq{EZstar0} and rearranging terms, provided that all the terms are finite, it remains to show that when $c_n<\infty$ then $m_n<\infty$, hence $c_k,m_k<\infty$ for $k\le n$. Recall \eq{zetan} which implies that $\zeta_n$ is stochastically bounded by the process $V_0+\sum_{i=1}^n V_i\prod_{j=1}^iU_j$ and thus $Z^*\le_{\text{st}} V_0+\sum_{i=1}^\infty V_i\prod_{j=1}^iU_j$, where in this case $V_i\sim\W_\tau$. One may apply Minkowski's inequality to argue that if $EV_0^n=EW_\tau^n<\infty$, then
    \begin{equation}
    \|V_0+\sum_{i=1}^n V_i\prod_{j=1}^iU_j\|_n\le \|V_0\|_n\sum_{i=0}^n\|U_1\|_n^i\le\frac{\|V_0\|}{1-\|U\|_n}
    \end{equation}
    where $\|Y\|=(E|Y|^n)^{ 1/n}$, so that by monotone convergence we also have that $E(V_0+\sum_{i=1}^\infty V_i\prod_{j=1}^iU_j)^n<\infty$ and thus $E(Z^*)^n<\infty$. To see that $EW_\tau^n<\infty$ whenever $c_n<\infty$ we recall that $X_\tau=W_\tau-L_\tau$ where $W_\tau,L_\tau$ are independent and $L_\tau\sim\exp(\alpha_\lambda)$ has moments of all orders. Since whenever $c_n<\infty$, $EX_t^n$ is a polynomial in $t$ of order at most $n$, which is easy to check by first computing factorial moments via $\log Ee^{-\alpha X_t}=\varphi(\alpha)t$, and since $\tau\sim\exp(\lambda)$ has finite moments of all orders, then $EX^n_\tau<\infty$ and therefore necessarily also $EW_\tau^n<\infty$, which completes the argument.
\end{proof}

It turns out that when $U\sim\text{Uniform}(0,1)$ the problem of identifying $Ee^{-\alpha Z^*}$ for all $\alpha\ge 0$ (in particular $Ef(\alpha_\lambda U)/\alpha_\lambda$), becomes tractable. For this case $EU^n=(n+1)^{-1}$.

\begin{theorem}
\label{thm2}
    Assume that $U\sim\text{Uniform}(0,1)$ and that $X$ is a finite mean spectrally positive L\'evy process, which is not a subordinator. Then, with
    \begin{equation}\label{eq:b}
    b=\left(\int_0^{\alpha_\lambda} e^{-\int_0^x\frac{\varphi(y)}{y(\lambda-\varphi(y))}dy}dx\right)^{-1},
    \end{equation}
    we have that
\begin{equation}\label{eq:f_alpha}
  f(\alpha)=Ee^{-\alpha Z^*}=\begin{cases}\frac{b}{1-\frac{\varphi(\alpha)}{\lambda}}\,\int_\alpha^{\alpha_\lambda}e^{-\int_\alpha^x\frac{\varphi(y)}{y(\lambda-\varphi(y))}dy}dx, &\alpha\in[0,\alpha_\lambda),\\ \\ 
  \frac{b}{\frac{1}{\alpha_\lambda}+\frac{\varphi'(\alpha_\lambda)}{\lambda}},&\alpha=\alpha_\lambda,\\ \\
  \frac{b}{\frac{\varphi(\alpha)}{\lambda}-1}\,\int_{\alpha_\lambda}^\alpha e^{-\int_x^\alpha\frac{\varphi(y)}{y(\varphi(y)-\lambda)}dy}dx,&\alpha\in(\alpha_\lambda,\infty).
  \end{cases}
\end{equation}
 Moreover, whenever $X_t=J_t-dt$ where $J$ is a nonzero pure jump subordinator and $d>0$ then $P(Z^*=0)=\frac{\lambda b}{2d}$. In all other cases $P(Z^*=0)=0$.
\end{theorem}

\medskip
\begin{proof}
    We first note that when $U\sim\text{Uniform}(0,1)$ then with $F(\alpha)=\int_0^\alpha f(x)dx$, $b=2F(\alpha_\lambda)/\alpha_\lambda^2$ and $a(\alpha)=\left(\alpha\left(1-\frac{\varphi(\alpha)}{\lambda}\right)\right)^{-1}$, \eq{falpha} becomes
    \begin{equation}\label{eq:diffeq}
    f(\alpha)=F'(\alpha)=\frac{F(\alpha)-\frac{\alpha^2}{2}b}{\alpha\left(1-\frac{\varphi(\alpha)}{\lambda}\right)}=a(\alpha)F(\alpha)-\frac{\alpha^2}{2}a(\alpha)b\,.
    \end{equation}
    From this it immediately follows via L'H\^opital's rule on the right hand side and the fact that $\varphi(\alpha_\lambda)=\lambda$, that
    \begin{equation}
    f(\alpha_\lambda)=\frac{f(\alpha_\lambda)-\alpha_\lambda b}{-\alpha_\lambda\frac{\varphi'(\alpha_\lambda)}{\lambda}}\,,
    \end{equation}
    so that the case $\alpha=\alpha_\lambda$ in \eq{f_alpha} is immediate.

    From \eq{diffeq} it now follows that, for $0<\alpha<\beta<\alpha_\lambda$ and for $\alpha_\lambda<\alpha<\beta$,
    \begin{equation}\label{eq:sol}
F(\beta)e^{-\int_\alpha^\beta a(y)dy}=F(\alpha)-b\int_\alpha^\beta \frac{x^2}{2}a(x)e^{-\int_\alpha^xa(y)dy}dx\,.
\end{equation}
Integration by parts gives
\begin{equation}
\int_\alpha^\beta \frac{x^2}{2}a(x)e^{-\int_\alpha^xa(y)dy}dx=\frac{\alpha^2}{2}-\frac{\beta^2}{2} e^{-\int_\alpha^\beta a(y)dy}+
\int_\alpha^\beta xe^{-\int_\alpha^x a(y)dy}dx\,.
\end{equation}
Since for a sufficiently small neighborhood of $\alpha_\lambda$ we have that 
\begin{equation}
\frac{\varphi(\alpha)-\lambda}{\alpha-\alpha_\lambda}\in (\varphi'(\alpha_\lambda)-\epsilon,\varphi'(\alpha_\lambda)+\epsilon)\,,
\end{equation}
where $0<\epsilon<\varphi'(\alpha_\lambda)$ and $\varphi'(\alpha_\lambda)$ is finite and positive,  then for every $\alpha\in (0,\alpha_\lambda)$ we necessarily have that $\int_\alpha^{\alpha_\lambda}a(x)dx=\infty$.
Therefore, letting $\beta\uparrow\alpha_\lambda$ implies that
\begin{equation}
F(\alpha)=b\left( \frac{\alpha^2}{2}+\int_\alpha^{\alpha_\lambda}xe^{-\int_\alpha^{x}a(y)dy}dx\right)\,.
\end{equation}
Differentiating both sides leads to
\begin{equation}\label{eq:f}
f(\alpha)=Ee^{-\alpha Z^*}=b a(\alpha)\int_\alpha^{\alpha_\lambda}xe^{-\int_\alpha^{x}a(y)dy}dx .
\end{equation}
We now recall that $a(\alpha)=(\alpha (1-\varphi(\alpha)/\lambda))^{-1}$ and note that for $\alpha<x<\alpha_\lambda$ we have that $x/\alpha=e^{\int_\alpha^x \frac{1}{y}dy}$, implying that
\begin{equation}
a(\alpha)xe^{-\int_\alpha^x a(y)dy}=\frac{1}{1-\frac{\varphi(\alpha)}{\lambda}}\, e^{\int_\alpha^x \left(\frac{1}{y}-a(y)\right)dy}=\frac{1}{1-\frac{\varphi(\alpha)}{\lambda}}
\,e^{-\int_\alpha^x \frac{\varphi(y)}{y(\lambda-\varphi(y))}dy}\,,
\end{equation}
so that together with \eq{f}, the case $\alpha\in [0,\alpha_\lambda)$ in \eq{f_alpha} follows.

Recalling, by assumption, that $\varphi'(0)=-EX_1$ is finite then $\int_0^x \frac{\varphi(y)}{y(\lambda-\varphi(y))}dy$ is convergent for any $x\in (0,\alpha_\lambda)$ so that, since $\varphi(\alpha)\to 0$ as $\alpha\downarrow 0$, we finally obtain \eq{b}.

For the case $\alpha>\alpha_\lambda$ we observe that in \eq{sol}
we can multiply both sides by $e^{\int_\alpha^\beta a(y)dy}$ and notice that for $y>\alpha_\lambda$ 
\begin{equation}
    a(y)=\frac{1}{y\left(1-\frac{\varphi(y)}{\lambda}\right)}=-\frac{1}{y\left(\frac{\varphi(y)}{\lambda}-1\right)}\,.
\end{equation}
The rest is a repetition of the proof for the case $\alpha\in[0,\alpha_\lambda)$, where eventually we let $\alpha\downarrow \alpha_\lambda$ and replace $\beta$ in the resulting expression by $\alpha$.

We finally turn to $P(Z^*=0)$. Note that 
\begin{equation}
   \int_{\alpha_\lambda}^\alpha e^{-\int_x^\alpha \frac{\varphi(y)}{y(\varphi(y)-\lambda)}dy}dx =
   \alpha\int_{\alpha_\lambda/\alpha}^1 e^{-\int_x^1 \frac{\varphi(\alpha y)}{y(\varphi(\alpha y)-\lambda)}dy}dx\,,
\end{equation}
and that by bounded convergence
\begin{equation}  \lim_{\alpha\to\infty}\int_{\alpha_\lambda/\alpha}^1 e^{-\int_x^1 \frac{\varphi(\alpha y)}{y(\varphi(\alpha y)-\lambda)}dy}dx=\int_0^1e^{-\int_x^1\frac{1}{y}dy}dx=\int_0^1xdx=\frac{1}{2}\,.
    \end{equation}
    Therefore,
    \begin{equation}\label{eq:zero}
    P(Z^*=0)=\frac{b}{2}\lim_{\alpha\to\infty}\frac{\alpha}{\frac{\varphi(\alpha)}{\lambda}-1}=\frac{\lambda b}{2}\,\lim_{\alpha\to\infty}\frac{\alpha}{\varphi(\alpha)}\,.
\end{equation}
If our L\'evy process is a subordinator minus a positive drift $d$ then
\begin{equation}
    \varphi(\alpha)=d\alpha-\int_{(0,\infty)}\left(1-e^{-\alpha x}\right)\nu(dx)=\alpha\left(d-\int_0^\infty e^{-\alpha x}\nu(x,\infty)dx\right)\,.
\end{equation}
In this case the mean is finite iff $\int_0^\infty \nu(x,\infty)dx=\int_{(0,\infty)}x\nu(dx)<\infty$, in which case (bounded convergence) $\int_0^\infty e^{-\alpha x}\nu(x,\infty)dx$ vanishes, whence, $\varphi(\alpha)/\alpha\to d$ as $\alpha\to\infty$.

In all other cases (recalling that $X$ is not a subordinator), that is, when there is a Brownian motion part or $\int_{(0,1]}x\nu(dx)=\infty$ (that is, an unbounded variation component) we have that $\varphi(\alpha)/\alpha\to \infty$ as $\alpha\to\infty$, hence the right hand side of \eq{zero} is zero.

\end{proof}

\begin{remark}\label{rem:g}{\rm

It is easy to check that for $\alpha\in[0,\alpha_\lambda)$ and with
\begin{equation}
    g(\alpha):=\int_0^{\alpha} e^{-\int_0^x \frac{\varphi(y)}{y(\lambda-\varphi(y))}dy}dx\,,
\end{equation}
we have the following representation (when $U\sim\text{Uniform}(0,1)$).
\begin{align}\label{eq:fal}
    f(\alpha)=Ee^{-\alpha Z^*}=\frac{1-\frac{g(\alpha)}{g(\alpha_\lambda)}}{g'(\alpha)\left(1-\frac{\varphi(\alpha)}{\lambda}\right)}.
\end{align}
Consequently, from \eq{wzutau} and \eq{LSTwtau} it follows that
\begin{equation}
    Ee^{-\alpha (Z^*U-L_\tau)^+}=\frac{1-\frac{g(\alpha)}{g(\alpha_\lambda)}}{g'(\alpha)\left(1-\frac{\alpha}{\alpha_\lambda}\right)}\,.
\end{equation}
We complete this remark by emphasizing that $b=1/g(\alpha_\lambda)$.
}
\end{remark}

\begin{remark}\label{rem:beta}{\rm
The proof of Theorem~\ref{thm2} can be easily modified to allow the more general case of $U\sim\text{Beta}(\theta,1)$. For $b$ and the cases $\alpha\in[0,\alpha_\lambda)$ and $\alpha\in(\alpha_\lambda,\infty)$, the resulting expressions are the same as in \eq{b} and \eq{f_alpha}, except  that the integrals in the exponents are multiplied by $\theta$. For the case $\alpha=\alpha_\lambda$, $1/\alpha_\lambda$ in the denominator is replaced by $\theta/\alpha_\lambda$. The expression $P(Z^*=0)=\frac{\lambda b}{2d}$ for the case of a subordinator with a negative drift is replace by $\frac{\lambda b}{(1+\theta)d}$. One may also apply the same ideas to handle a mixture of a $\text{Beta}(\theta,1)$ distribution and an atom at zero (=clearing). For the sake of brevity we refrain from working out these exercises here.
}    
\end{remark}

\begin{remark}
{\rm
Regarding Remark~\ref{rem:beta}, we mention in passing that if we consider an on/off process with $\exp(\lambda)$ distributed on times and $\exp(\eta)$ distributed off times, where during on times the process behaves like a reflected L\'evy process and during off times the process is a deterministic linear dam with release rate $r>0$ (decreases at rate $rx$ where $x$ is a current state) then this process restricted to on times is the process $Z$ and hence has the steady-state distribution of $Z^*$ with $U\sim\text{Beta}(r/\eta,1)$. Furthermore, by PASTA, the steady-state distribution of the process restricted to off times is the distribution of $Z^*U$, the LST of which may be deduced from \eq{falpha}. Therefore the steady-state distribution of the on/off process is a $\left(\frac{\eta}{\lambda+\eta},\frac{\lambda}{\lambda+\eta}\right)$ mixture of the distributions of $Z^*$ and $Z^*U$. Thus, our analysis covers this case as well.
}
\end{remark}

\section{Examples and tail behavior}
\label{secexample}
In this section we work out the LST of  $Z^*$, and hence of $Z$, for the two cases of reflected Brownian motion with collapse (Subsection~\ref{sec6.1}) and
the $M/M/1$ queue with collapse (Subsection~\ref{sec6.2}). In both cases, collapses occur according to a Poisson process with rate $\lambda$, and instantaneously remove a Uniform$(0,1)$ distributed fraction of the workload present. In these two examples we will use the notation $g$ and the expression for $f$ from Remark~\ref{rem:g}. We conclude in Subsection~\ref{sec6.3} with a discussion about the tail behavior of $Z^*$ for the compound Poisson case with a negative drift.

\subsection{Reflected Brownian motion with collapse}
\label{sec6.1}
The Laplace exponent of Brownian motion with  drift $c$ is given by (cf.\ (\ref{Laplaceexp}))
\begin{equation}
\varphi(\alpha) = -c\alpha + \frac{\sigma^2}{2} \alpha^2 .
\end{equation}
In order to obtain the LST $f(\alpha)$ of the steady-state level of the reflected Brownian motion with collapses,
we need to determine $g(\alpha)$ and its derivative, which feature in (\ref{eq:fal}).
First write
\begin{equation}
   - \frac{\varphi(y)}{y(\lambda - \varphi(y))} ~=~ \frac{c - \frac{\sigma^2}{2} y}{\lambda +cy - \frac{\sigma^2}{2} y^2} ~=~ 
    \frac{y-\frac{2c}{\sigma^2}}{(y-y_1)(y-y_2)} ,
    \label{quot}
\end{equation}
where 
\begin{equation}
    y_{1,2} ~=~ \frac{c}{\sigma^2} \pm \sqrt{\frac{c^2}{\sigma^4} + \frac{2\lambda}{\sigma^2}} .
    \end{equation}
    It is easily verified that $y_1>0$ and $y_2<0$, with $y_1y_2 = -\frac{2\lambda}{\sigma^2}$ and $y_1+y_2 = \frac{2c}{\sigma^2}$; and of course $y_1 = \alpha_\lambda$, the unique positive zero of $\lambda - \varphi(\alpha)$.

    Using a partial fraction expansion, we can rewrite (\ref{quot}) as follows:
    \begin{equation}
    - \frac{\varphi(y)}{y(\lambda-\varphi(y))} ~=~ \frac{D_1}{y-y_1} + \frac{D_2}{y-y_2} ,
    \end{equation}
    where 
    \begin{equation}
        D_1 = \frac{-y_2}{y_1-y_2} > 0, ~~~~ D_2 = \frac{y_1}{y_1-y_2} > 0.
        \label{D1D2}
    \end{equation}
    It is now easily verified that
    \begin{equation}
        g(\alpha) ~=~\int_0^\alpha \left(\frac{x-y_1}{-y_1}\right)^{D_1} \left(\frac{x-y_2}{-y_2}\right)^{D_2} dx ,
    \end{equation}
    and hence
    \begin{equation}
        \frac{g(\alpha)}{g(\alpha_\lambda)} ~=~\frac{g(\alpha)}{g(y_1)}  ~=~
        \frac{\int_0^\alpha (x-y_1)^{D_1} (x-y_2)^{D_2} dx}{\int_0^{y_1} (x-y_1)^{D_1} (x-y_2)^{D_2} dx} ,
    \end{equation}
    while
    \begin{equation}
        g'(\alpha) ~=~ \left(\frac{\alpha -y_1}{-y_1}\right)^{D_1} \left(\frac{\alpha -y_2}{-y_2}\right)^{D_2} .
    \end{equation}
    Substitution in (\ref{eq:fal}) shows that the LST of the steady-state level of our reflected Brownian motion with collapse is, for $\alpha \in [0,y_1)$, given by 
    \begin{align}
    f(\alpha) &=
    \frac{y_1y_2}{(\alpha-y_1)(\alpha-y_2)}
    \left(\frac{\alpha -y_1}{-y_1}\right)^{-D_1} \left(\frac{\alpha -y_2}{-y_2}\right)^{-D_2} 
    \label{falBM} \nonumber \\ \nonumber\\
    &\qquad \times \left(1 - \frac{\int_0^\alpha (x-y_1)^{D_1} (x-y_2)^{D_2} dx}{\int_0^{y_1} (x-y_1)^{D_1} (x-y_2)^{D_2} dx} \right) 
    \\ \nonumber\\
    &= \left(\frac{y_1}{y_1 - \alpha}\right)^{1+D_1} \left(\frac{-y_2}{-y_2 + \alpha}\right)^{1+D_2} 
    \left(1 -\frac{\int_0^\alpha (x-y_1)^{D_1} (x-y_2)^{D_2} dx}{\int_0^{y_1} (x-y_1)^{D_1} (x-y_2)^{D_2} dx} \right) .
    \nonumber
    \end{align}
    Notice that the terms before the large brackets, in the last line of (\ref{falBM}), represent the LST of the difference of two Gamma distributed random variables.
    \begin{remark}
    \label{remarkbeta}
    It should be observed that $g(\alpha)$ can be expressed in incomplete beta functions.
    The incomplete beta function $B(z;a_1,a_2)$ is defined as (e.g., p.~263 of \cite{AS1970})
    \begin{equation}
        B(z;a_1,a_2) ~:=~\int_0^z t^{a_1-1} (1-t)^{a_2-1} dt ;
        \label{Betaf}
    \end{equation}
    it can also be expressed as a hypergeometric function via the identity
    \begin{equation}
        B(z;a_1,a_2) ~=~ \frac{z^{a_1}}{a_1}  ~ {}_2F_1(a_1,1-a_2;a_1+1;z) .
        \end{equation}
        The transformation $t:= \frac{x-y_2}{y_1-y_2}$ gives:
        \begin{eqnarray}
            && \int_0^\alpha (x-y_1)^{D_1} (x-y_2)^{D_2} dx ~=~ 
            (-1)^{D_1} (y_1-y_2)^{D_1+D_2+1} 
 ~\int_{\frac{-y_2}{y_1-y_2}}^{\frac{\alpha-y_2}{y_1-y_2}} t^{D_2}(1-t)^{D_1} dt
 \nonumber\\
 &=& (-1)^{D_1} (y_1-y_2)^2 \left[B\left(\frac{\alpha-y_2}{y_1-y_2};1+D_2,1+D_1\right) - B\left(\frac{-y_2}{y_1-y_2};1+D_2,1+D_1\right) \right].
 \nonumber\\
        \end{eqnarray}
        \end{remark}

\subsection{The M/M/1 queue with collapse}
\label{sec6.2}
A second example in which one can explicitly determine the LST of the steady-state reflected spectrally-positive L\'evy process 
with collapse concerns the workload in an $M/M/1$ queue with collapse. Denote the arrival rate of customers by $\gamma$ and the service rate by $\mu$; the server works at speed $d$ work units per time unit. Collapses again occur according to a Poisson process with rate $\lambda$, and remove (and maintain) a Uniform$(0,1)$ distributed part of the workload present.

The Laplace exponent of the workload process in an $M/M/1$ queue is given by
\begin{equation}
    \varphi(\alpha) ~=~ d \alpha - \gamma \frac{\alpha}{\mu+\alpha} .
    \label{varphimm1}
    \end{equation}
Very similar to the previous subsection, we can write:
\begin{equation}
    - \frac{\varphi(y)}{y(\lambda - \varphi(y))} ~=~
\frac{y+\mu - \frac{\gamma}{d}}{(y-z_1)(y-z_2)},
\label{quot2}
\end{equation}
where
\begin{equation}
    z_{1,2} ~=~\frac{1}{2} \left[ \frac{\lambda+\gamma}{d} - \mu \pm \sqrt{\left(\frac{\lambda+\gamma}{d}-\mu\right)^2 + 4 \frac{\lambda\mu}{d} } \,\right].
\end{equation}
It is easily verified that $z_1>0$ and $z_2<0$, with $z_1z_2 = - \frac{\lambda \mu}{d}$ and $z_1+z_2=\frac{\lambda+\gamma}{d}-\mu$; and of course $z_1=\alpha_\lambda$, the unique positive zero of $\lambda - \varphi(\alpha)$.

Using a partial fraction expansion, we can rewrite (\ref{quot2}) as follows:
\begin{equation}
    - \frac{\varphi(y)}{y(\lambda - \varphi(y))} ~=~
\frac{F_1}{y-z_1} + \frac{F_2}{y-z_2},
\label{quot3}
\end{equation}
where
\begin{equation}
    F_1 = \frac{\frac{\lambda}{d}-z_2}{z_1-z_2} > 0, ~~~~ F_2 = \frac{z_1 - \frac{\lambda}{d}}{z_1-z_2} >0.
\end{equation}
It is now easily verified that
    \begin{equation}
        g(\alpha) ~=~\int_0^\alpha \left(\frac{x-z_1}{-z_1}\right)^{F_1} \left(\frac{x-z_2}{-z_2}\right)^{F_2} dx ,
    \end{equation}
    and hence
    \begin{equation}
        \frac{g(\alpha)}{g(\alpha_\lambda)} ~=~\frac{g(\alpha)}{g(z_1)}  ~=~
        \frac{\int_0^\alpha (x-z_1)^{F_1} (x-z_2)^{F_2} dx}{\int_0^{z_1} (x-z_1)^{F_1} (x-z_2)^{F_2} dx} ,
    \end{equation}
    while
    \begin{equation}
        g'(\alpha) ~=~ (\frac{\alpha -z_1}{-z_1})^{F_1} (\frac{\alpha -z_2}{-z_2})^{F_2} .
    \end{equation}
    Substitution in (\ref{eq:fal}) shows that the LST of the steady-state workload in the $M/M/1$ queue with collapse is, for $\alpha \in [0,z_1)$,  given by 
   \begin{align}
    f(\alpha) &=
    \frac{z_1z_2}{(\alpha-z_1)(\alpha-z_2)}
    \left(\frac{\alpha -z_1}{-z_1}\right)^{-F_1} \left(\frac{\alpha -z_2}{-z_2}\right)^{-F_2} 
    \label{falMM1} \nonumber \\ \nonumber \\
    &\qquad\times ~\left(1 - \frac{\int_0^\alpha (x-z_1)^{F_1} (x-z_2)^{F_2} dx}{\int_0^{z_1} (x-z_1)^{F_1} (x-z_2)^{F_2} dx} \right) 
    \\ \nonumber\\
    &= \left(\frac{z_1}{z_1 - \alpha}\right)^{1+F_1} \left(\frac{-z_2}{-z_2 + \alpha}\right)^{1+F_2} 
    \left(1 -\frac{\int_0^\alpha (x-z_1)^{F_1} (x-z_2)^{F_2} dx}{\int_0^{z_1} (x-z_1)^{F_1} (x-z_2)^{F_2} dx} \right) .
    \nonumber
    \end{align} 
    Just as in Remark~\ref{remarkbeta}, one can express this LST in terms of incomplete beta functions. Furthermore, the terms before the large brackets in the last line of (\ref{falMM1}) is the LST of the difference of two Gamma distributed random variables.
\begin{remark}
    \label{remarkmm1}{\rm
    We briefly indicate another approach that basically leads to the same first-order inhomogeneous differential equation as \eq{diffeq}.
The steady-state workload distribution has an atom $p_0$ at zero and conjecture that it is absolutely continuous on $(0,\infty)$ with density $h(\cdot)$, so that $P(Z^*\le t)=p_0+\int_0^t h(s)\,ds$. Level crossing arguments lead to the following equation:
\begin{equation}
d h(x)+ \lambda \int_x^\infty h(y)\,\frac{x}{y}\,dy= \gamma \int_0^x e^{-\mu(x-y)}dP(Z^* <y) \,.
\label{LCT}
\end{equation}
Here the first term in the left hand side corresponds to downcrossings of level $x$ due to serving, and the second term in the left hand side reflects downward jumps due to a uniform collapse with rate $\lambda$. The term in the right hand side corresponds to a jump through level $x$ due to a customer arrival -- either a jump from zero or from some level $y \in (0,x)$. Multiplying both sides of (\ref{LCT}) by $e^{-\alpha x}$ and integrating from zero to infinity yields, with $f(\alpha) = Ee^{-\alpha Z^*} = p_0 + \int_0^\infty e^{-\alpha x} h(x) dx$:
\begin{equation} 
d (f(\alpha) -p_0)  + \lambda \int_{y=0}^\infty \frac{h(y)}{y} \int_{x=0}^y x e^{-\alpha x} dx dy = \frac{\gamma}{\mu+\alpha} f(\alpha) .
\end{equation}
The second term in the left hand side can be rewritten as
\begin{align}
    \lambda \int_{y=0}^\infty h(y)\, \frac{1 - e^{-\alpha y} -\alpha y e^{-\alpha y}}{\alpha^2 y} dy &= \frac{\lambda}{\alpha^2} \int_{u=0}^\alpha (f(u)-p_0) du - \lambda\, \frac{f(\alpha)-p_0}{\alpha}\nonumber\\ &= \frac{\lambda}{\alpha^2} F(\alpha) - \frac{\lambda}{\alpha} f(\alpha) ,
\end{align}
where $F(\alpha) = \int_0^\alpha f(x) dx$, cf.\ above \eq{diffeq}, and hence we have
\begin{equation}
    F'(\alpha) \left(\frac{\lambda}{\alpha} -d + \frac{\gamma}{\mu+\alpha}\right) = \frac{\lambda}{\alpha^2} F(\alpha)  - d p_0 .
    \end{equation}
Dividing both sides by
 $\frac{\lambda}{\alpha} -d + \frac{\gamma}{\mu+\alpha}$ results in \eq{diffeq}, with $b=2dp_0/\lambda$ -- which is in agreement with Theorem~\ref{thm2}. 
 }
 \end{remark}       

\subsection{Tail behavior for compound Poisson with negative drift}
\label{sec6.3}
Consider the special case in which the spectrally positive L\'evy process consists of a compound Poisson process minus a positive drift.
We denote the arrival rate of the Poisson jumps by $\gamma$, and a generic jump size by $B$, with distribution $B(\cdot)$ and LST $\beta(\cdot)$; the drift is denoted by $d$.
The resulting model is an $M/G/1$ queue with arrival rate $\gamma$, generic service requirement $B$, service speed $d$, and uniformly distributed collapses occurring according to a Poisson process with rate $\lambda$. $Z^*$ represents the steady-state workload of this queue; its LST is for $\alpha \in [0,\alpha_\lambda)$ given by (\ref{eq:fal}), with $\varphi(\alpha) = d\alpha - \gamma(1-\beta(\alpha))$.

We now demonstrate how (\ref{eq:fal}) can be used to determine the tail behavior\break $P(Z^* > t)$ for large $t$, in the case of regularly varying jump sizes of index\break $-\delta \in (-2,-1)$.
In this case, the mean service requirement is finite but the variance is infinite. The Bingham-Doney lemma (cf.\ Theorem 8.1.6 of \cite{BGT}) states that (i) and (ii) below are equivalent:
\begin{description}
\item{(i)} $P(B>t)  \approx  \frac{-1}{\Gamma(1-\delta)} t^{-\delta}\, l(t)$, as $t \rightarrow \infty$,
\item{(ii)} $\beta(\alpha) -1 +EB \alpha ~\approx ~\alpha^{\delta}\, l(1/\alpha)$, as $\alpha \downarrow 0$.
\end{description}
Here $l(\cdot)$ is a slowly varying function at infinity.
If (i) holds, and hence (ii), then we have, for $\alpha \downarrow 0$,
\begin{equation}
    \frac{\lambda}{\lambda - \varphi(\alpha)} 
 - 1 - (d-\gamma EB)\alpha ~\approx ~\frac{\gamma}{\lambda} \alpha^{\delta} l(1/\alpha) ,
 \end{equation}
 \begin{equation}
 \int_0^\alpha \frac{\varphi(y)}{y(\lambda - \varphi(y))} dy ~-~\frac{d-\gamma EB}{\lambda} \alpha ~\approx ~ \frac{\gamma}{\delta \lambda} \alpha^{\delta} l(1/\alpha),
 \end{equation}
 and hence, recalling Remark~\ref{rem:g}, for $\alpha \downarrow 0$,
 \begin{equation}
     g'(\alpha) -1 +  \frac{d-\gamma EB}{\lambda} \alpha ~\approx ~ - \frac{\gamma}{\delta \lambda} \alpha^{\delta},
 \end{equation}
 \begin{equation}
     g(\alpha) - \alpha ~\approx ~- \frac{d - \gamma EB}{2 \lambda} \alpha^2 .
 \end{equation}
Combining these asymptotic expansions with (\ref{eq:fal}) shows that
\begin{equation}
Ee^{-\alpha Z^*} - 1 + (\frac{1}{g(\alpha_\lambda)} -2 \frac{d - \gamma EB}{\lambda}) \alpha ~\approx ~ \frac{\gamma}{\lambda} \frac{\delta+1}{\delta} \alpha^{\delta} l(1/\alpha), ~~~ \alpha \downarrow 0.
    \label{RV3a}
    \end{equation}
Now observe, using \eq{EZstar1}, that $\frac{1}{g(\alpha_\lambda)} - 2 \frac{d-\gamma EB}{\lambda}$ equals $EZ^*$. Indeed, $\frac{1}{g(\alpha_\lambda)} = b$, and on the other hand the term $\frac{1}{1-EU} \frac{Ef(\alpha_\lambda U)}{\alpha_\lambda}$ in \eq{EZstar1} can be rewritten in the Uniform$(0,1)$ case as $\frac{2}{\alpha_\lambda^2} \int_0^{\alpha_\lambda} f(y) dy = \frac{2}{\alpha_\lambda^2} F(\alpha_\lambda)$, which equals $b$ (see above Equation \eq{diffeq}). 
Hence we have
    \begin{equation}
    Ee^{-\alpha Z^*} - 1 + EZ^* \alpha ~\approx ~ \frac{\gamma}{\lambda} \frac{\delta+1}{\delta} \alpha^{\delta} l(1/\alpha), ~~~ \alpha \downarrow 0.
    \label{RV3}
    \end{equation}
    Applying the Bingham-Doney lemma in the reverse direction reveals that
    \begin{equation}
        P(Z^*>t) ~\approx ~ \frac{-1}{\Gamma(1-\delta)} \frac{\gamma}{\lambda} \frac{\delta+1}{\delta} t^{-\delta} l(t), ~~~t \rightarrow \infty .
    \end{equation}
    In an ordinary $M/G/1$ queue with regularly varying service requirements, the workload tail is one degree heavier than the tail of the service requirement,
    whereas in this $M/G/1$ queue with uniform collapses, $Z^*$ is regularly varying of index $-\delta$, just as $B$: $P(Z^*>t) \approx  \frac{\gamma}{\lambda} \frac{\delta+1}{\delta} P(B>t)$.
 In the $M/G/1$ shot-noise queue (an $M/G/1$ queue with service speed proportional to the workload) the same phenomenon was observed in \cite{BM}. This is not so surprising, as uniform collapses can also be seen as proportional (shot-noise) decrements during an exponentially distributed amount of time.

\end{document}